\documentclass[12pt]{amsart}
\usepackage{amssymb}
\usepackage[dvips]{graphics}
\ifx\mathrm\undefined\let\mathrm\rm\fi
\ifx\mathbf\undefined\let\mathbf\bf\fi
\ifx\mathfrak\undefined\let\mathfrak\frak\fi
\ifx\mathcal\undefined\let\mathcal\cal\fi
\ifx\mathbb\undefined\let\mathbb\Bbb\fi
\ifx\emph\undefined\let\emph\it\fi
\newcommand{\g}{{{\mathfrak g}\,}}

\newcommand{\C}{{\mathbb C}}
\newcommand{\Q}{{\mathbb Q}}

\newcommand{\Ref}[1]{{(\ref{#1})}}
\newcommand{\be}{\begin{displaymath}}
\newcommand{\ee}{\end{displaymath}}
\newcommand{\bea}{\begin{eqnarray*}}
\newcommand{\eea}{\end{eqnarray*}}

\newcommand{\w}{{W_{sl_2}}}
\newcommand{\T}{{\mathcal{T}}}
\newcommand{\U}{{\mathcal{U}}}

 at9.98pt

\newcommand{\dontprint}[1]{\relax}

\newtheorem%
{thm}{Theorem}[section]
\newtheorem%
{proposition}[thm]{Proposition}
\newtheorem%
{lemma}[thm]{Lemma}
\newtheorem%
{lemmadef}[thm]{Lemma-Definition}
\newtheorem%
{corollary}[thm]{Corollary}
\newtheorem%
{conjecture}[thm]{Conjecture}

\title[The $sl_2$ approximation of the Kontsevich integral
of the unknot]{A remark on the $sl_2$ approximation of the Kontsevich integral
of the unknot}
\author[S. Tyurina and A. Varchenko]
{Svetlana Tyurina${}^{*}$ and Alexander Varchenko${}^{**,1}$}
\thanks{${}^1$Supported in part by NSF grant  DMS-9801582}
\date{June 2001}
\begin{document}
%
%
\def\bcirc
{    \bezier{40}(-1,0)(-1,0.577)(-0.5,0.866)
     \bezier{8}(-0.5,0.866)(0,1.155)(0.5,0.866)
     \bezier{40}(0.5,0.866)(1,0.577)(1,0)
     \bezier{8}(0.5,-0.866)(1,-0.577)(1,0)
     \bezier{40}(-0.5,-0.866)(0,-1.155)(0.5,-0.866)
     \bezier{8}(-1,0)(-1,-0.577)(-0.5,-0.866)
}
\def\pone{\put(-0.866,0.5){\circle*{0.15}}}
\def\ptwo{\put(0.866,0.5){\circle*{0.15}}}
\def\pthree{\put(0,-1){\circle*{0.15}}}
\def\pfour{\put(-0.25,-0.966){\circle*{0.15}}}
\def\pfive{\put(0.25,-0.966){\circle*{0.15}}}
\def\psix{\put(0.707,0.707){\circle*{0.15}}}
\def\pseven{\put(0.966,0.25){\circle*{0.15}}}
\def\dftrone{\begin{picture}(0,0) \bcirc \pone\ptwo\pfour\pfive
     \bezier{80}(-0.866,0.5)(0,0)(0.25,-0.966)
     \bezier{80}(0.866,0.5)(0,0)(-0.25,-0.966) \end{picture} }
\def\dftrtwo{\begin{picture}(0,0) \bcirc \pone\ptwo\pfour\pfive
     \bezier{80}(-0.866,0.5)(-0.28,-0.12)(-0.25,-0.966)
     \bezier{80}(0.866,0.5)(0.28,-0.12)(0.25,-0.966) \end{picture} }
\def\dftrthree{\begin{picture}(0,0) \bcirc \pone\pthree\psix\pseven
     \bezier{80}(-0.866,0.5)(-0.04,0.3)(0.707,0.707)
     \bezier{80}(0,-1)(0.24,-0.19)(0.966,0.25) \end{picture} }
\def\dftrfour{\begin{picture}(0,0) \bcirc \pone\pthree\psix\pseven
     \bezier{80}(-0.866,0.5)(0,0)(0.966,0.25)
     \bezier{80}(0,-1)(0,0)(0.707,0.707) \end{picture} }
%
%
%
%
\def\trpax
{
    \put(-0.866,0.5){\circle*{0.15}}
    \put(-0.866,-0.5){\circle*{0.15}}
    \put(0.866,0.5){\circle*{0.15}}
    \put(0.866,-0.5){\circle*{0.15}}
    \put(0,1){\circle*{0.15}}
    \put(0,-1){\circle*{0.15}}
    \put(0,-1){\line(0,1){2}}
    \bezier{80}(-0.866,0.5)(0,0.2)(0.866,0.5)
    \bezier{80}(-0.866,-0.5)(0,-0.2)(0.866,-0.5)
}
\def\trp
{   \begin{picture}(40,40)
    \put(0,0){\circle{2}}
    \trpax
    \end{picture}
}
\def\fctw
{  \bezier{25}(-0.26,0.97)(0,1.035)(0.26,0.97)
   \bezier{25}(0.26,0.97)(0.52,0.9)(0.71,0.71)
   \bezier{25}(0.71,0.71)(0.9,0.52)(0.97,0.26)
   \bezier{25}(0.97,-0.26)(0.9,-0.52)(0.71,-0.71)
   \bezier{25}(0.71,-0.71)(0.52,-0.9)(0.26,-0.97)
   \bezier{25}(0.26,-0.97)(0,-1.035)(-0.26,-0.97)
   \bezier{25}(-0.71,-0.71)(-0.9,-0.52)(-0.97,-0.26)
   \bezier{25}(-0.97,-0.26)(-1.035,0)(-0.97,0.26)
   \bezier{25}(-0.97,0.26)(-0.9,0.52)(-0.71,0.71)
   \bezier{25}(-0.71,0.71)(-0.52,0.9)(-0.26,0.97)
}
%
%
\def\horch{\put(-1,0){\circle*{0.2}}\put(1,0){\circle*{0.2}}
    \put(-1,0){\line(1,0){2}}}
\def\lch{\put(-0.707,0.707){\circle*{0.2}}\put(-0.707,-0.707){\circle*{0.2}}
    \bezier{30}(-0.707,-0.707)(-0.3,0)(-0.707,0.707)}
\def\rch{\put(0.707,-0.707){\circle*{0.2}}\put(0.707,0.707){\circle*{0.2}}
    \bezier{30}(0.707,-0.707)(0.3,0)(0.707,0.707)}
\def\done{}
\def\dtwo{\lch}
\def\dthree{\rch}
\def\dfour{\lch\rch}
\def\dfive{\horch}
\def\dsix{\horch\lch}
\def\dseven{\horch\rch}
\def\deight{\horch\lch\rch}
\def\diag#1{\unitlength=10pt \put(0,0){\circle{2}} #1}
\def\cdl#1{\raisebox{0pt}[16pt][8pt]{ 
    \begin{picture}(20,10)(-8,-3.5) \diag{#1} \end{picture}} }
\def\cdsl#1{\mbox{\begin{picture}(26,10)(-12,-3.5) \diag{#1} \end{picture}} }
\def\wcdl#1#2{ #1 \Bigl( \raisebox{0pt}[16pt][8pt]{ 
    \begin{picture}(16,10)(-6,-3.5) \diag{#2} \end{picture}} \Bigr) }
\def\wcdsl#1#2{ #1 \Bigl( \mbox{ 
    \begin{picture}(16,10)(-6,-3.5) \diag{#2} \end{picture}} \Bigr) }
\def\verch{
  \put(0,1){\circle*{0.2}}\put(0,-1){\circle*{0.2}}\put(0,-1){\line(0,1){2}}}
\def\dtrkr{\verch\dkrc}
\def\chlnpv{\put(-0.707,-0.707){\circle*{0.2}}
            \put(0,1){\circle*{0.2}}
            \put(-0.707,-0.707){\line(2,5){0.7}} }
\def\chlvpn{\put(-0.707,0.707){\circle*{0.2}}
            \put(0,-1){\circle*{0.2}}
            \put(-0.707,0.707){\line(2,-5){0.7}} }
\def\dtrtr{\chlnpv \chlvpn \rch}
\def\dtrfo{\rch \lch \verch}
\def\luch{\put(-0.707,0.707){\circle*{0.2}}
          \put(0,1){\circle*{0.2}}
          \bezier{20}(-0.707,0.707)(-0.21,0.45)(0,1)}
\def\ldch{\put(-0.707,-0.707){\circle*{0.2}}
          \put(0,-1){\circle*{0.2}}
          \bezier{20}(-0.707,-0.707)(-0.21,-0.45)(0,-1)}
\def\dtrfi{\luch \ldch \rch}
%
%
%
\def\onesix{\put(1,0){\circle*{0.2}} \put(-0.707,-0.707){\circle*{0.2}}
             \put(1,0){\line(-5,-2){1.8}} }
\def\threig{\put(0,1){\circle*{0.2}} \put(0.707,-0.707){\circle*{0.2}}
             \put(0.707,-0.707){\line(-2,5){0.7}} }
\def\twofiv{\put(-1,0){\circle*{0.2}} \put(0.707,0.707){\circle*{0.2}}
             \put(0.707,0.707){\line(-5,-2){1.8}} }
\def\twosev{\put(0,-1){\circle*{0.2}} \put(0.707,0.707){\circle*{0.2}}
             \put(0.707,0.707){\line(-2,-5){0.7}} }
\def\fiveig{\put(-1,0){\circle*{0.2}} \put(0.707,-0.707){\circle*{0.2}}
             \put(0.707,-0.707){\line(-5,2){1.8}} }
\def\onethr{\put(1,0){\circle*{0.2}}\put(0,1){\circle*{0.2}}
    \bezier{30}(1,0)(0.23,0.23)(0,1)}
\def\onesev{\put(1,0){\circle*{0.2}}\put(0,-1){\circle*{0.2}}
    \bezier{30}(1,0)(0.23,-0.23)(0,-1)}
\def\fivsev{\put(-1,0){\circle*{0.2}}\put(0,-1){\circle*{0.2}}
    \bezier{30}(-1,0)(-0.23,-0.23)(0,-1)}
\def\foufiv{\put(-1,0){\circle*{0.2}}\put(-0.707,0.707){\circle*{0.2}}
          \bezier{20}(-1,0)(-0.45,0.21)(-0.707,0.707)}
\def\oneeig{\put(1,0){\circle*{0.2}}\put(0.707,-0.707){\circle*{0.2}}
          \bezier{20}(1,0)(0.45,-0.21)(0.707,-0.707)}
\def\twothr{\put(0,1){\circle*{0.2}}\put(0.707,0.707){\circle*{0.2}}
          \bezier{20}(0,1)(0.21,0.45)(0.707,0.707)}
%
%
%
\def\dfoone{\onesix \twofiv \threig \chlvpn}
\def\dfotwo{\onethr \twosev \fiveig \lch}
\def\dfothree{\onesev \fiveig \lch \twothr}
\def\dfofour{\onethr \rch \fivsev \lch}
\def\dfofive{\oneeig \twothr \lch \fivsev}
\def\dfosix{\oneeig \twothr \foufiv \ldch}
\def\dfoseven{\horch \twosix \verch \foueig}
\def\dfoeight{\onesix \twofiv \verch \foueig}
\def\dfonine{\horch \rch \chlnpv \chlvpn}
\def\dfoten{\horch \rch \verch \lch}
%
%
%
\def\bepi#1{\makebox[36pt]{\unitlength=20pt 
            \begin{picture}(1.8,1.1)(-0.98,-0.2)  #1
            \end{picture}} }
\def\sctw#1#2#3#4{
   \bezier{25}(-0.26,0.97)(0,1.035)(0.26,0.97)
   \bezier{4}(0.26,0.97)(0.52,0.9)(0.71,0.71)
   \bezier{#1}(0.71,0.71)(0.9,0.52)(0.97,0.26)
   \bezier{4}(0.97,0.26)(1.035,0)(0.97,-0.26)
   \bezier{#2}(0.97,-0.26)(0.9,-0.52)(0.71,-0.71)
   \bezier{4}(0.71,-0.71)(0.52,-0.9)(0.26,-0.97)
   \bezier{25}(0.26,-0.97)(0,-1.035)(-0.26,-0.97)
   \bezier{4}(-0.26,-0.97)(-0.52,-0.9)(-0.71,-0.71)
   \bezier{#3}(-0.71,-0.71)(-0.9,-0.52)(-0.97,-0.26)
   \bezier{4}(-0.97,-0.26)(-1.035,0)(-0.97,0.26)
   \bezier{#4}(-0.97,0.26)(-0.9,0.52)(-0.71,0.71)
   \bezier{4}(-0.71,0.71)(-0.52,0.9)(-0.26,0.97)
}
%
%
%
\def\uch{\put(-0.707,0.707){\circle*{0.2}}\put(0.707,0.707){\circle*{0.2}}
    \bezier{30}(-0.707,0.707)(0,0.3)(0.707,0.707) }
\def\dch{\put(0.707,-0.707){\circle*{0.2}}\put(-0.707,-0.707){\circle*{0.2}}
    \bezier{30}(-0.707,-0.707)(0,-0.3)(0.707,-0.707) }
\def\dparh{\uch \dch}
\def\dtrpar{\verch \uch \dch}
%
%
%
\def\twosix{\put(0.707,0.707){\circle*{0.2}}\put(-0.707,-0.707){\circle*{0.2}}
    \put(-0.707,-0.707){\line(1,1){1.414}} }
\def\foueig{\put(-0.707,0.707){\circle*{0.2}}\put(0.707,-0.707){\circle*{0.2}}
    \put(-0.707,0.707){\line(1,-1){1.414}} }
\def\dkrc{\twosix \foueig}

\maketitle
\centerline{\it ${}^{*}$MPIM, Bonn,}
\centerline{\it Germany}
\centerline{tyurina@mpim-bonn.mpg.de}
\medskip
\centerline{\it ${}^{**}$Department of Mathematics,
University of North Carolina at Chapel Hill,}
\centerline{\it Chapel Hill, NC 27599-3250, USA}
\centerline{av@math.unc.edu}
\newcommand{\s}{\sigma}
\begin{abstract}
${}$

The Kontsevich integral of a knot $K$ is a sum 
$I(K)=1+\sum_{n=1}^\infty h^n\sum_{D\in A_n}a_D D$ 
over all chord diagrams with suitable coefficients. 
Here $A_n$ is the space of chord diagrams with $n$ chords.
A simple explicit formula for coefficients $a_D$ is not 
known even for the unknot.
Let $E_1, E_2, \ldots$ be elements of $A=\oplus_{n}A_n.$
Say that a sum $I'(K)=1+\sum_{n=1}^\infty h^n E_n$ is an $sl_2$ 
approximation of the Kontsevich integral if the values 
of the $sl_2$ weight system $\w$ on both sums are equal,
$\w (I(K))=\w (I'(K))$.

For any natural $n$ fix points $a_1,...,a_{2n}$ on a circle. For any
permutation  $\sigma \in S^{2n}$ of $2n$ elements 
define the chord diagram $D(\sigma)$ with $n$ chords as the diagram
with chords formed by pairs  $a_{\s(2i-1)}$ and $a_{\s(2i)},$
$i=1,\ldots , n$. We show that
$$
1+\sum_{n=1}^\infty \frac{h^{2n}}{2^n(2n)!(2n+1)!}\sum_{\sigma \in S^{2n}}D(\sigma)
$$
is an $sl_2$ approximation of the Kontsevich integral of the unknot.
\end{abstract}

\section{ $sl_2$ approximations }

The Kontsevich integral of a knot $K$ is a sum $I(K)=1+\sum_{n=1}^\infty h^n\sum_{D\in A_n}a_D D$ 
over all chord diagrams with suitable coefficients \cite{K}. 
Here $A_n$ is the space of chord diagrams with $n$ chords, $h$ is formal parameter.
A simple explicit formula for coefficients $a_D$ is not 
known even for the unknot.
Let $E_1, E_2, \ldots$ be elements of $A=\oplus_{n}A_n.$
Say that a sum $I'(K)=1+\sum_{n=1}^\infty h^n E_n$ is an $sl_2$ 
approximation of the Kontsevich integral if the values 
of the $sl_2$ weight system $\w$ on both sums are equal,
$\w (I(K))=\w (I'(K))$.

For any natural $n$ fix points $a_1,...,a_{2n}$ on a circle. For any
permutation  $\sigma \in S^{2n}$ of $2n$ elements 
define the chord diagram $D(\sigma)$ with $n$ chords as the diagram
with chords formed by pairs $a_{\s(2i-1)}$ and $a_{\s(2i)},$
$i=1,\ldots , n$. Introduce an element 
$\Sigma_n=\sum_{\s \in S^{2n}}D(\s)$ in $A_n$.
\begin{thm}\label{main}
The  sum
\begin{equation}\label{main}
1+\sum_{n=1}^\infty \frac{h^{2n}}{2^n(2n)!(2n+1)!}\Sigma_n 
\notag
\end{equation}
is an $sl_2$ approximation of the Kontsevich integral of the unknot.
\end{thm}
The theorem is proved in Section 6.

We thank S.Chmutov for useful discussions. 

\section{Three algebras}

The Kontsevich integral takes values in the graded completion 
of the chord diagram algebra $A=\oplus_n A_n$. 
The $\Q$ vector space  $A$
is generated by the usual chord diagrams 
modulo the four-term relation

%
%
\begin{center}
\begin{picture}(300,50)(-20,-25)
  \setlength{\unitlength}{20pt}
  \put(0,0){\dftrone}
  \put(1.9,0){\makebox(0,0){$-$}}
  \put(3.7,0){\dftrtwo}
  \put(5.6,0){\makebox(0,0){$+$}}
  \put(7.4,0){\dftrthree}
  \put(9.3,0){\makebox(0,0){$-$}}
  \put(11.1,0){\dftrfour}
  \put(13.1,0){\makebox(0,0){$=$}}
  \put(14.1,0){\makebox(0,0){$0\,.$}}
\end{picture}
\end{center}

The product of chord diagrams is their connected  sum,  see \cite{Bar95}. 

We also consider the algebra $\T$ of trivalent diagrams.
A trivalent diagram is a connected 
graph with only trivalent vertices and a distinguished 
oriented circle, such that at each vertex, which does not 
lie on the circle, one of two possible cyclic 
orderings of the three edges meeting at this vertex is chosen. 
The $\Q$ vector space $\T$ is
generated by trivalent diagrams modulo the STU relation, 

\begin{picture}(400,40)
\put(50,20){\circle*{3}}
\put(50,10){\circle*{3}}
\put(50,10){\line(0,1){10}}
\put(50,20){\line(1,1){10}}
\put(50,20){\line(-1,1){10}}
\put(90,15){$=$}
\put(40,10){\vector(1,0){20}}
\put(40,10.51){\vector(1,0){20}}
\put(40,9.51){\vector(1,0){20}}
\put(55,0){${}_S$}

\put(150,10){\circle*{3}}
\put(160,10){\circle*{3}}
\put(150,10){\line(-1,2){10}}
\put(160,10){\line(1,2){10}}
\put(175,15){$-$}
\put(140,10){\vector(1,0){30}}
\put(140,10.51){\vector(1,0){30}}
\put(140,9.51){\vector(1,0){30}}
\put(160,0){${}_T$}
\put(200,10){\circle*{3}}
\put(210,10){\circle*{3}}
\put(200,10){\line(1,2){10}}
\put(210,10){\line(-1,2){10}}
\put(190,10){\vector(1,0){30}}
\put(190,10.51){\vector(1,0){30}}
\put(190,9.51){\vector(1,0){30}}
\put(210,0){${}_U$}
\put(230,10){$.$}
\end{picture}
\newline
These three trivalent diagrams are 
identical outside the corresponding fragment
on the picture. Pieces of the circle are pictured by thick lines. 
The product of trivalent diagrams is defined as their connected sum
with respect to the distinguished circles.

\medskip
The following relations follow from the STU relation,
%
%



\begin{center}
   \begin{picture}(300,40)(0,-20)
       \unitlength=20pt
   \put(0,0.5){\makebox(0,0){(AS):}}
   \put(3,0){
      \begin{picture}(0,0)
      \put(0,0){\circle*{0.15}}
      \put(0,0){\line(0,-1){1}}
      \put(0,0){\line(-1,1){1}}
      \put(0,0){\line(1,1){1}}
      \end{picture}
              }
   \put(5,0){\makebox(0,0){$=$}}
   \put(6,0.05){\makebox(0,0){$-$}}
   \put(7,0){
      \begin{picture}(0,0)
      \put(0,0){\circle*{0.15}}
      \put(0,0){\line(0,-1){1}}
      \put(-1,1){\line(2,-1){1}}
      \bezier{60}(0,0.5)(0.3,0.3)(0,0)
      \put(1,1){\line(-2,-1){1}}
      \bezier{60}(0,0.5)(-0.3,0.3)(0,0)
      \end{picture}
    \put(0.5,0){\makebox(0,0){$,$}}          
              }
      \end{picture}
\end{center}
%
%
\def\angl{
      \put(0,0){\circle*{0.15}}
      \put(0,0){\line(-1,0){1}}
      \put(0,0){\line(3,1){1.5}}
      \put(0,0){\line(3,-1){1.5}}
         }
\def\ihx#1#2{
\begin{center}
   \begin{picture}(300,40)(0,-20)
       \unitlength=20pt
   \put(0,0.45){\makebox(0,0){(#1):}}
   \put(3,0){
      \begin{picture}(0,0)
      #2
      \put(-0.5,0){\circle*{0.15}}
      \put(-0.5,0){\line(0,-1){1}}
      \end{picture}
              }
   \put(6,0){\makebox(0,0){$=$}}
   \put(8,0){
      \begin{picture}(0,0)
      #2
      \put(0.75,0.25){\circle*{0.15}}
      \put(0.75,0.25){\line(0,-1){1.25}}
      \end{picture}
              }
   \put(10.2,0){\makebox(0,0){$+$}}
   \put(11.7,0){
      \begin{picture}(0,0)
      #2
      \put(0.75,-0.25){\circle*{0.15}}
      \put(0.75,-0.25){\line(0,-1){0.75}}
      \end{picture}
              }
   \put(14,0){\makebox(0,0){$.$}}
   \end{picture}
\end{center}
}
%
%
\ihx{IHX}{\angl}

Applying the STU relation one can express a 
given trivalent diagram as a linear combination of chord diagrams. 
This gives a natural mapping $\T \to A$ which is
 an isomorphism  of algebras, see \cite{Bar95}.

The third algebra is the algebra $\U$ of uni-trivalent diagrams. 
A uni-trivalent diagram is a graph whose vertices either 
univalent or trivalent, at each trivalent vertex a cyclic 
ordering of the three edges is chosen.
The vertices of valency 1 of a uni-trivalent diagram are called ``legs'' of
the diagram. 
We consider the $\Q$ vector space $\U$
generated by uni-trivalent graphs
modulo the AS and IHX relations.
The product in the algebra $\U$ is the disjoint union $\sqcup$. 
There is a natural isomorphism $\U \to \T$ as vector spaces, but not as algebras. The isomorphism
maps every uni-trivalent
diagram to the average of all possible ways of 
placing its univalent vertices along the circle, see \cite{Bar95}.

Introduce uni-trivalent graphs $w_{2n}$, the ``wheels with $2n$ legs'',

\begin{picture}(400,40)
\put(0,15){$w_2=$}
\put(60,15){,}
\put(40,20){\circle{15}}
\put(45,25){\line(1,1){10}}
\put(35,15){\line(-1,-1){10}}

\put(160,15){$, \qquad \ldots$}


\put(100,15){$w_4=$}
\put(140,20){\circle{15}}
\put(145,25){\line(1,1){10}}
\put(135,15){\line(-1,-1){10}}
\put(145,15){\line(1,-1){10}}
\put(135,25){\line(-1,1){10}}
\end{picture}

Under the isomorphism $\U\to\T$, for instance, we have

\begin{picture}(200,40)
\put(40,15){$\mapsto \frac{1}{2}\,(\,{}\,$}
\put(160,15){$\,).$}
\put(110,15){$+$}
\put(90,15){\circle{15}}
\put(85,0){\line(0,1){10}}
\put(95,0){\line(0,1){10}}
\put(90,15){\circle{30}}
\put(140,15){\circle{15}}
\put(136,0){\line(1,1){10}}
\put(144,0){\line(-1,1){10}}
\put(140,15){\circle{30}}
\put(25,20){\circle{15}}
\put(30,25){\line(1,1){10}}
\put(20,15){\line(-1,-1){10}}
\end{picture}

\section{The $sl_2$-weight system}
A weight system on $A$ (resp. on $\T,\,\U$) with values in a vector
space is a linear homomorphism of $A$  
(resp. of $\T,\,\U$) to the vector space.
The composition of the Kontsevich integral of a knot with the linear
homomorphism defines a knot invariant of the knot.
Here we recall the construction of the weight system associated to a
Lie algebra $\g$ with an $ad$-invariant nondegenerate bilinear form $\Phi$. 

Let ${a_1,\ldots,a_m}$ and ${b_1,\ldots,b_m}$
be two dual bases of $\g$: $\Phi(a_i,b_j) = \delta_{i,j}$.
Fix a chord diagram $D$ and a base point on its circle 
which is different from the endpoints of the chords. 
Label each chord by a number $i$
such that $1\leq i \leq m$. Attach to one endpoint of the chord labelled
by $i$ the element $a_i$ and to another endpoint the element $b_i$. 
For example, for the chord diagram 
%
%
$D = \mbox{
\begin{picture}(30,18)(-8,-2.5)
    \unitlength=10pt
    \put(0,0){\circle{2}}
    \put(-0.2,1.05){\raisebox{-7pt}{*}}
    \put(-0.707,0.707){\circle*{0.2}}
    \put(0.707,-0.707){\circle*{0.2}}
    \put(-0.707,0.707){\line(1,-1){1.414}}
    \bezier{5}(0.65,-0.45)(1.1,-0.45)(1.5,-0.45)
    \put(1.7,-0.4){\makebox(0,0){\scriptsize{$i$}}}
    \put(0.707,0.707){\circle*{0.2}}
    \put(-0.707,-0.707){\circle*{0.2}}
    \put(-0.707,-0.707){\line(1,1){1.414}}
    \bezier{5}(0.65,0.45)(1.1,0.45)(1.5,0.45)
    \put(1.7,0.45){\makebox(0,0){\scriptsize{$j$}}}
\end{picture}}$ 
with the base point  \raisebox{-4pt}{*} and labels $i$ and $j$ we 
have
%
%
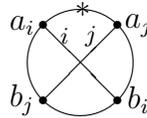
\begin{center}
\begin{picture}(40,30)(-20,-15)
    \unitlength=20pt
    \put(0,0){\circle{2}}
    \put(-0.1,1.045){\raisebox{-7pt}{*}}
    \put(-0.707,0.707){\circle*{0.15}}
    \put(0.707,-0.707){\circle*{0.15}}
    \put(-0.707,0.707){\line(1,-1){1.414}}
    \put(-0.3,0.5){\makebox(0,0){\scriptsize{$i$}}}
    \put(0.707,0.707){\circle*{0.15}}
    \put(-0.707,-0.707){\circle*{0.15}}
    \put(-0.707,-0.707){\line(1,1){1.414}}
    \put(0.2,0.5){\makebox(0,0){\scriptsize{$j$}}}
    \put(-1.1,0.707){\makebox(0,0){$a_i$}}
    \put(-1.1,-0.707){\makebox(0,0){$b_j$}}
    \put(1.1,0.707){\makebox(0,0){$a_j$}}
    \put(1.1,-0.707){\makebox(0,0){$b_i$}}
\end{picture}
\end{center}
Walk around the circle starting from the base point in the
direction of the orientation of the circle and write 
in one word the elements associated to the endpoints .
 The constructed
word is an element of the
universal enveloping algebra $U(\g)$. Define
$W(D)$ to be the sum of such words where the sum is 
over all labels of the chords. In
the example,
$$ W(D) = \sum_{i,j} a_i b_j b_i a_j \in U(\g) . $$ 
In our pictures we always assume  that the circle
is oriented counterclockwise.

The element $W(D)$ does not depend on the base point,
does not depend on the choice of dual bases in $\g$,
belongs to the center $Z(\g)$ of the universal 
enveloping algebra, 
satisfies the four-term relation. The mapping
$W : A\to Z(\g)$  is an algebra homomorphism, 
$W(D_1 \cdot D_2) = W(D_1) \cdot W(D_2)$, see \cite{K}.

{\bf Example.} $\wcdl{W}{\dfive} = c$ is the quadratic Casimir element of 
$Z(\g)$ associated to the chosen invariant form.

For $\g=sl_2$, we have $Z(sl_2) \cong \C[c]$, and for a chord diagram $D$ 
with $n$ chords, 
$$\w(D) = c^n + \lambda_1 c^{n-1} + \lambda_2 c^{n-2} + \ldots + 
\lambda_{n-1} c.$$

 We choose $Tr$ as an $ad$-invariant form on $sl_2$ where
$Tr$ is the trace of matrices in 
the standard two dimensional representation of $sl_2.$
A recurrent formula for $\w$ is constructed in \cite{CV}.


%
%
\def\bepi#1{\makebox[36pt]{\unitlength=20pt 
            \begin{picture}(1.8,1.1)(-0.98,-0.2)  #1
            \end{picture}} }
\def\sctw#1#2#3#4{
   \bezier{25}(-0.26,0.97)(0,1.035)(0.26,0.97)
   \bezier{4}(0.26,0.97)(0.52,0.9)(0.71,0.71)
   \bezier{#1}(0.71,0.71)(0.9,0.52)(0.97,0.26)
   \bezier{4}(0.97,0.26)(1.035,0)(0.97,-0.26)
   \bezier{#2}(0.97,-0.26)(0.9,-0.52)(0.71,-0.71)
   \bezier{4}(0.71,-0.71)(0.52,-0.9)(0.26,-0.97)
   \bezier{25}(0.26,-0.97)(0,-1.035)(-0.26,-0.97)
   \bezier{4}(-0.26,-0.97)(-0.52,-0.9)(-0.71,-0.71)
   \bezier{#3}(-0.71,-0.71)(-0.9,-0.52)(-0.97,-0.26)
   \bezier{4}(-0.97,-0.26)(-1.035,0)(-0.97,0.26)
   \bezier{#4}(-0.97,0.26)(-0.9,0.52)(-0.71,0.71)
   \bezier{4}(-0.71,0.71)(-0.52,0.9)(-0.26,0.97)
}
%
%
\def\cirtw
{  \sctw{25}{25}{25}{25}
   \put(0.866,0.5){\circle*{0.15}}
   \put(0.866,-0.5){\circle*{0.15}}
   \put(-0.866,-0.5){\circle*{0.15}}
   \put(-0.866,0.5){\circle*{0.15}}
   \put(1.25,0.5){\makebox(0,0){$e'_i$}}
   \put(1.25,-0.5){\makebox(0,0){$e'_j$}}
   \put(-1.2,-0.5){\makebox(0,0){$e_j$}}
   \put(-1.2,0.5){\makebox(0,0){$e_i$}}
}

\begin{thm}

%
%
\def\slchone{\put(0,1){\circle*{0.15}}\put(0,-1){\circle*{0.15}}
             \put(0,-1){\line(0,1){2}} }
\def\slchtwo{\put(-0.174,0.985){\circle*{0.15}}
             \put(0.866,0.5){\circle*{0.15}}
             \bezier{70}(-0.174,0.985)(0.21,0.45)(0.866,0.5) }
\def\slchthree{\put(-0.174,-0.985){\circle*{0.15}}
               \put(0.866,-0.5){\circle*{0.15}}
               \bezier{70}(-0.174,-0.985)(0.21,-0.45)(0.866,-0.5) }
\def\slchfour{\put(0.174,0.985){\circle*{0.15}}
              \put(0.866,0.5){\circle*{0.15}}
              \bezier{60}(0.174,0.985)(0.36,0.48)(0.866,0.5) }
\def\slchfive{\put(0.174,-0.985){\circle*{0.15}}
              \put(0.866,-0.5){\circle*{0.15}}
              \bezier{60}(0.174,-0.985)(0.36,-0.48)(0.866,-0.5) }
\def\slchsix{\put(-0.174,0.985){\circle*{0.15}}
             \put(-0.174,-0.985){\circle*{0.15}}
             \bezier{100}(-0.174,0.985)(0.3,0)(-0.174,-0.985) }
\def\slchseven{\put(-0.174,-0.985){\circle*{0.15}}
               \put(0.866,0.5){\circle*{0.15}}
               \bezier{90}(-0.174,-0.985)(0,0)(0.866,0.5) }
\def\slcheight{\put(-0.174,0.985){\circle*{0.15}}
               \put(0.866,-0.5){\circle*{0.15}}
               \bezier{90}(-0.174,0.985)(0,0)(0.866,-0.5) }
\def\slchnine{\put(-0.174,0.985){\circle*{0.15}}
              \put(-0.866,0.5){\circle*{0.15}}
              \bezier{60}(-0.174,0.985)(-0.36,0.48)(-0.866,0.5) }
\def\slchten{\put(0.174,0.985){\circle*{0.15}}
             \put(-0.866,0.5){\circle*{0.15}}
             \bezier{90}(0.174,0.985)(-0.21,0.45)(-0.866,0.5) }
\def\slchel{\put(-0.14,0.985){\circle*{0.15}}
            \put(0.174,-0.985){\circle*{0.15}}
            \put(0.2,-0.985){\line(-1,6){0.325}} }
\def\slchtw{\put(0.174,-0.985){\circle*{0.15}}
            \put(-0.866,0.5){\circle*{0.15}}
            \bezier{90}(0.174,-0.985)(0,0)(-0.866,0.5) }
\def\slchtth{\put(0.866,0.5){\circle*{0.15}}
             \put(0.866,-0.5){\circle*{0.15}}
             \bezier{60}(0.866,0.5)(0.5,0)(0.866,-0.5) }
%
%
\def\achone{\put(0.174,0.985){\circle*{0.15}}
            \put(0.866,-0.5){\circle*{0.15}}
            \bezier{80}(0.174,0.985)(0.36,0.16)(0.866,-0.5) }
\def\achtwo{\put(0.174,-0.985){\circle*{0.15}}
            \put(0.866,0.5){\circle*{0.15}}
            \bezier{80}(0.174,-0.985)(0.36,-0.16)(0.866,0.5) }
\def\achthree{\put(-0.174,0.985){\circle*{0.15}}
            \put(-0.866,-0.5){\circle*{0.15}}
            \bezier{80}(-0.174,0.985)(-0.36,0.16)(-0.866,-0.5) }
\def\achfour{\put(-0.174,-0.985){\circle*{0.15}}
            \put(-0.866,0.5){\circle*{0.15}}
            \bezier{80}(-0.174,-0.985)(-0.36,-0.16)(-0.866,0.5) }
\def\achfive{\put(-0.174,-0.985){\circle*{0.15}}
              \put(-0.866,-0.5){\circle*{0.15}}
              \bezier{60}(-0.174,-0.985)(-0.36,-0.48)(-0.866,-0.5) }
\def\achsix{\put(0.174,-0.985){\circle*{0.15}}
            \put(-0.866,-0.5){\circle*{0.15}}
            \bezier{80}(0.174,-0.985)(-0.21,-0.45)(-0.866,-0.5) }
\def\achseven{\put(0.174,0.985){\circle*{0.15}}
             \put(0.174,-0.985){\circle*{0.15}}
             \bezier{100}(0.174,0.985)(-0.3,0)(0.174,-0.985) }
\def\acheight{\put(-0.866,0.5){\circle*{0.15}}
             \put(-0.866,-0.5){\circle*{0.15}}
             \bezier{60}(-0.866,0.5)(-0.5,0)(-0.866,-0.5) }
\def\achnine{\put(0.174,0.985){\circle*{0.15}}
             \put(-0.866,-0.5){\circle*{0.15}}
             \bezier{90}(0.174,0.985)(0,0)(-0.866,-0.5) }
%
%
\def\dslone{\bepi{\sctw{25}{25}{4}{4} \slchone \slchtwo \slchthree}}
\def\dsltwo{\bepi{\sctw{25}{25}{4}{4} \slchone \slchfour \slchthree}}
\def\dslthree{\bepi{\sctw{25}{25}{4}{4} \slchone \slchtwo \slchfive}}
\def\dslfour{\bepi{\sctw{25}{25}{4}{4} \slchone \slchfour \slchfive}}
\def\dslfive{\bepi{\sctw{25}{25}{4}{4} \slchsix \slchtth}}
\def\dslsix{\bepi{\sctw{25}{25}{4}{4} \slchseven \slcheight}}
\def\dslseven{\bepi{\sctw{4}{25}{4}{25} \slchone \slchnine \slchthree}}
\def\dsleight{\bepi{\sctw{4}{25}{4}{25} \slchone \slchten \slchthree}}
\def\dslnine{\bepi{\sctw{4}{25}{4}{25} \slchone \slchnine \slchfive}}
\def\dslten{\bepi{\sctw{4}{25}{4}{25} \slchone \slchten \slchfive}}
\def\dslel{\bepi{\sctw{4}{25}{4}{25} \slchel \put(-0.866,0.5){\circle*{0.15}}
    \put(0.866,-0.5){\circle*{0.15}} \put(0.866,-0.5){\line(-5,3){1.72}} }}
\def\dsltw{\bepi{\sctw{4}{25}{4}{25} \slchtw \slcheight}}
%
%
\def\aslone{\bepi{\sctw{25}{25}{4}{4} \slchone \slchseven \slcheight}}
\def\asltwo{\bepi{\sctw{25}{25}{4}{4} \slchone \slchseven \achone}}
\def\aslthree{\bepi{\sctw{25}{25}{4}{4} \slchone \achtwo \slcheight}}
\def\aslfour{\bepi{\sctw{25}{25}{4}{4} \slchone \achone \achtwo}}
\def\aslfive{\bepi{\sctw{25}{25}{4}{4} \slchtwo \slchthree}}
\def\asltse{\bepi{\sctw{25}{4}{25}{4} \slchone \achsix \slchtwo}}
\def\asltei{\bepi{\sctw{25}{4}{25}{4} \slchone \achsix \slchfour}}
\def\asltni{\bepi{\sctw{25}{4}{25}{4} \slchone \achfive \slchtwo}}
\def\asltte{\bepi{\sctw{25}{4}{25}{4} \slchone \achfive \slchfour}}
\def\asltto{\bepi{\sctw{25}{4}{25}{4} \achnine \slchseven}}
\def\aslttt{\bepi{\sctw{25}{4}{25}{4}\slchel\put(-0.866,-0.5){\circle*{0.15}}
   \put(0.866,0.5){\circle*{0.15}}\put(0.866,0.5){\line(-5,-3){1.72}} }}
%
%
\def\di#1{\raisebox{0pt}[20pt][23pt]{#1}}
\def\wdi#1{W\biggl( {\raisebox{0pt}[20pt][30pt]{#1}} \biggr)}
\def\spwdi#1{W\biggl( \mbox{#1} \biggr)}
\def\slrel#1#2#3#4#5#6#7{
\begin{eqnarray*}
\wdi{#1}-\wdi{#2}-\wdi{#3}+\wdi{#4} = \makebox[30pt]{}&  \\
\hfill = 2 \wdi{#5} - 2 \wdi{#6} #7 &
\end{eqnarray*} }

Let $W=\w$ be the weight system associated to 
$sl_2$ and the $ad$-invariant form $Tr$. Then
\slrel{\dslone}{\dsltwo}{\dslthree}{\dslfour}{\dslfive}{\dslsix};
\slrel{\aslone}{\asltwo}{\aslthree}{\aslfour}{\dslfive}{\aslfive};
\slrel{\dsleight}{\dslseven}{\dslten}{\dslnine}{\dsltw}{\dslel};
\slrel{\asltse}{\asltei}{\asltni}{\asltte}{\asltto}{\aslttt}.
\end{thm}
This theorem allows one to
compute $W(D)$ since the two chord diagrams of
the right hand side have one chord less than the diagrams of the left hand
side, and the last three 
diagrams of the left hand side are simpler than the
first one since they have less intersections between their chords.

\medskip

The theorem indicates six-term elements of the kernel of the $sl_2$ weight
system. The subspace $I$ of the algebra $A$ generated by the six-term
elements forms an ideal. The quotient algebra $A/I$ is generated by
two elements  $\cdl{\dfive}$ and  $\cdl{\dkrc}$. The ideal generated by
the six term elements and the element
 $\cdl{\dkrc} + 2\cdl{\dfive} -  \cdl{\dfive} \cdot \cdl{\dfive}$
is the kernel of the $sl_2$ weight system.

The linear isomorphisms 
$\U \to \T$ and $\T \to A$ induce weight systems 
$\w:\T \to Z(sl_2),$ $\w:\U \to Z(sl_2).$

\begin{thm}{\rm \cite{CV}}
The weight system $\w$ satisfies the following three term relation

\begin{picture}(90,30)
\put(0,5){$\w($}
\put(30,0){\line(1,1){10}}
\put(50,10){\line(1,1){10}}
\put(60,0){\line(-1,1){10}}
\put(40,10){\line(-1,1){10}}
\put(40,10){\line(1,0){10}}
\put(65,5){$)=$} 
\end {picture}
\begin{picture}(90,30)
\put(-10,5){$2\w($}
\put(30,-2.5){\line(1,1){10}}
\put(50,12.5){\line(1,1){10}}
\put(60,-2.5){\line(-1,1){10}}
\put(40,12.5){\line(-1,1){10}}
\put(40,12.5){\line(1,0){10}}
\put(40,7.5){\line(1,0){10}}
\put(65,5){$)-$} 
\end {picture}
\begin{picture}(90,30)
\put(-10,5){$2\w($}
\put(30,0){\line(3,2){30}}
\put(60,0){\line(-3,2){30}}
\put(65,5){$) $} 
\end {picture}

\noindent 
for any uni-trivalent diagrams differed only by the pictured fragments. 
\end{thm}
%
\def\bubl
{   \put(-0.97,0){\circle*{0.2}}
    \put(0.97,0){\circle*{0.2}}
    \bezier{15}(-0.97,0)(-0.9,0.26)(-0.71,0.45)
    \bezier{15}(-0.71,0.45)(-0.52,0.64)(-0.26,0.71)
    \bezier{15}(-0.26,0.71)(0,0.775)(0.26,0.71)
    \bezier{15}(0.26,0.71)(0.52,0.64)(0.71,0.45)
    \bezier{15}(0.71,0.45)(0.9,0.26)(0.97,0)
    \bezier{15}(-0.97,0)(-0.9,-0.26)(-0.71,-0.45)
    \bezier{15}(-0.71,-0.45)(-0.52,-0.64)(-0.26,-0.71)
    \bezier{15}(-0.26,-0.71)(0,-0.775)(0.26,-0.71)
    \bezier{15}(0.26,-0.71)(0.52,-0.64)(0.71,-0.45)
    \bezier{15}(0.71,-0.45)(0.9,-0.26)(0.97,0)
}
\begin{corollary}\label{c1}\cite{CV}
$$\w\left( \mbox{
\begin{picture}(28,15)(-13.5,-2.5)
    \unitlength=10pt
    \bubl
    \put(-2,0){\line(1,0){1.03}}
    \put(1,0){\line(1,0){1.03}}
 \end{picture} } 
\right) = 
4\w\left( \mbox{
\begin{picture}(28,15)(-13.5,-2.5)
    \unitlength=10pt
    \put(-2,0){\line(1,0){4}}
 \end{picture} } 
\right).
$$
\end{corollary}

\section{Bernoulli numbers and Bernoulli polynomials}
The modified Bernoulli numbers are defined by the series 
$$
\sum_{n=0}^{\infty} b_{2n} x^{2n} = \frac {1}{2} \ln
\frac {e^{x/2}-e^{-x/2}}{x/2}\,.
$$

The Bernoulli polynomials $B_n(x)$ are defined by the series
$$
\frac{ze^{zx}}{e^z-1}=\sum_{n=0}^{\infty}B_n(x)\frac{z^n}{n!}\,.
$$
The polynomial  $B_n(x)$ has degree $n$. Its top coefficient equals 1. 

Following \cite{LM} introduce the shifted Bernoulli polynomials
$q_n(x)$ by the condition
$$
q_n(\frac{x^2-1}{2})=\frac{2}{(2n+1)!}\frac{B_{2n+1}(\frac{1+x}{2})}{x}\,.
$$
The polynomial  $q_n(x)$ has degree $n$. Its top coefficient equals
$\frac{1}{2^n(2n+1)!}$ . 

\begin{thm}\label{wheel} For any natural $n$, we have 
$$
\w(w_{2n})\,= \,2^{2n+1}\,(2n+1)!\,q_n(c)\,,
\qquad
\sum_{\sigma\in S^{2n}} \,\w (D(\sigma)) \,= \,2^n \,(2n)!\,(2n+1)!\, q_n(c)\,.
$$
\end{thm}
The theorem is proved in Section 6.

\section{The Kontsevich integral of the unknot}
A formula for the logarithm of the Kontsevich integral $I$ 
of the unknot in terms of wheels is given in \cite{Bar97, T}.

\begin{thm}\label{roz}
\begin{equation}\label{}
I \,= \, 1 + \exp \,(\sum_{n=1}^{\infty} b_{2n} h^{2n} w_{2n} )\,
 =\, 1 +
(\sum_{n=1}^{\infty} b_{2n}h^{2n} w_{2n}) 
+\frac {1}{2} (\sum_{n=1}^{\infty} b_{2n}h^{2n} w_{2n})^2 + \ldots 
\notag
\end{equation}
\end{thm}

The  value of the $sl_2$ weight
system on the Kontsevich integral of the unknot is  calculated in \cite{LM}.
\begin{thm}\label{lm}
\begin{equation}\label{}
\w(I)=\sum_{n=0}^\infty q_n(c)h^{2n}. 
\notag
\end{equation}
\end{thm}

\section{Proofs}
\begin{lemma}\label{t1}
For $n_1+n_2+\ldots+n_k=n$, we have  
$$\w(w_{2n_1}\sqcup w_{2n_2} \sqcup \ldots \sqcup w_{2n_k})=
\frac{2^{n_1+n_2+\cdots+n_k+k}}{(2n)!}
\w(\Sigma_n).$$ 
\end{lemma}
{\bf Proof of the lemma.} We prove the lemma 
for $k=1$, general case is similar.
The three term relation applied to a vertex of $w_{2n}$ gives
$$
\w(w_{2n})= 2\w(\,\{ \, {}\,|\,{}\, \}\, \sqcup w_{2n-2}) - 2\w(t_{2n-2}) 
$$
where $\{\,{}\,|\,{}\,\}$ is the uni-trivalent graph with one edge and two
univalent vertices and \newline
$
t_{2n-2}=
\begin{picture}(25,10)
\put(5,-4){\line(0,1){10}}
\put(10,-4){\line(0,1){10}}
\put(12,-4){${}_{\cdots}$}
\put(21,-4){\line(0,1){10}}
\put(0,6){\line(1,0){25}}
\put(0,-10){\hbox to 20pt{\upbracefill}}
\put(-5,-17){${}_{(2n-2)\, legs}$}
\end{picture} 
$\,.
\newline $ {}$
\newline
Application of  the three term relation to the first two vertices of
$t_{2n-2}$ gives
$$\w(t_{2n-2})\,=\,2\,\w(\{\,|\,\} \sqcup t_{2n-4})\,-\,2\,\w(\{|\} \sqcup
t_{2n-4})\,=\,0\,.
$$
Then Corollary \ref{c1} implies
$
\w(w_{2n})\,=\,2^{n+1}\,\w(\sqcup_{n}\{\,|\,\})
$
where 
$\sqcup_{n}\{\,|\,\}$ is the diagram with $n$ edges and $2n$ univalent vertices. 
Glueing legs of $\sqcup_{n}\{\,|\,\}$
to the circle in all possible ways and dividing the sum
of the resulting chord diagrams by $(2n)!$ gives the lemma for $k=1$,
\begin{equation}\label{1}
\w(w_{2n})=\frac{2^{n+1}}{(2n)!}\sum_{\sigma\in S^{2n}}\w(D_n(\sigma))\,.
\end{equation}
\hfill $\square$

{\bf Proof of Theorems \ref{main} and \ref{wheel}.} Theorems \ref{roz}, \ref{lm}
and Lemma \ref{t1} imply that $\w (\Sigma_n)\,=\, \text{const}\,
q_n(c)$. Both sides are polynomials in $c$ of degree $n$. Comparing
the coefficients of $c^n$ one gets Theorem \ref{main} and the second
 equality of Theorem \ref{wheel}. Equation \Ref{1} implies the first
equality of Theorem \ref{wheel}.
\hfill $\square$

\end{document}